\documentclass[10pt]{article}
\usepackage{graphicx}

\usepackage[utf8]{inputenc}
\usepackage[T1]{fontenc} %Required for inserting images
\usepackage{amsmath}
\usepackage{amsbsy}
\usepackage{amsthm}
\usepackage{amsfonts}
\usepackage{amssymb}
\usepackage{mathtools}
\usepackage{hyperref}
\usepackage{geometry}
\usepackage{subcaption}
\usepackage{csquotes} %Adjust the layout options.
\usepackage[english]{babel} 
\usepackage{fancyhdr}
\usepackage{cleveref}
\usepackage{relsize}
\usepackage{caption}
\usepackage{siunitx}
\usepackage{float}
\usepackage[multiple]{footmisc}
\usepackage{newunicodechar}
\newunicodechar{ℓ}{\ell}
\usepackage{booktabs}
\usepackage{xcolor}
\usepackage[style=authoryear]{biblatex}
\definecolor{citecolor}{RGB}{0, 0, 255}
\usepackage{hyperref}
\hypersetup{colorlinks,allcolors=black}
\usepackage{todonotes}
\title{\textbf{Stochastic ordering of extreme order statistics in Archimax copula}}
\author{Sarikul Islam\footnote{Corresponding author's E-mail id: \href{sarikul_phd_math@kgpian.iitkgp.ac.in}{sarikul\_phd\_math@kgpian.iitkgp.ac.in}} \, and\, Nitin Gupta\,\footnote{\href{nitin.gupta@maths.iitkgp.ac.in}{nitin.gupta@maths.iitkgp.ac.in}}\\  \\
 \textsuperscript{\footnotemark[1]\,\,\,\footnotemark[2]} Department of Mathematics, Indian Institute of Technology Kharagpur \,\\ Kharagpur, 721302, India.} 
 \date{}
 
\begin{document} \maketitle  \section*{Abstract} An extension of Archimax copula class in more than two random variables ( Multivariate ) was introduced in \parencite{jagr2011generalization} for describing dependency structures among random variables in higher dimension, and some properties of Archimax copula were explored in \parencite{charpentier2014multivariate}. In this article, some results for stochastic ordering of extreme order statistics in \parencite{li2015ordering} are generalized and proved in Archimax copula. Stochastic ordering of sample extremes for PHR models is generalized and proved in Archimax copula. Examples with graphical illustrations are also presented. 
\subsection*{Keywords } Stochastic ordering,  Majorization, Archimax copula, PHR models.
\section{Introduction}
A important class of copulas is Archimedean copulas, that can model the dependency structure among random variables using a single generator function. They have been used in various fields such as finance, insurance, hydrology etc. However, one of the limitations of Archimedean copulas is that they cannot capture the tail dependency of random variables, which is the probability that extreme events occurs simultaneously. To overcome this drawback, some extensions of Archimedean copulas have been proposed, such as the extreme-value copulas, Archimax copulas, meta copulas etc. The archimax copulas, introduced in \parencite{jagr2011generalization}, is a generalization of the Archimedean copulas that allows for different levels of dependence in the medium and the tails of the random variables and can be estimated by nonparametric methods, such as kernel density estimation or maximum likelihood estimation. Archimax copulas are determined by two function parameters: first is the stable tail dependence function $\ell$, and second is the Archimedean generator $\phi$. The intermediate instances as well as the asymptotic dependence and independence of random variables can be captured by the Archimax copulas. Furthermore, the copula parameters may be easily estimated and simulated because of its straightforward and tractable form. Order statistics, on the other hand, are employed in actuarial science, industrial engineering, survival analysis, financial engineering, and dependability theory. There are significant similarities between order statistics and $k$ out of $n$ systems. If at least $k$ of the $n$ components of a system work, then the system is $k$ out of $n$ with order statistics $X_{{n-k+1}:n}$ as its lifetime. Stochastic comparisons of order statistics is an important topic in dependability theory. For reference see \parencite{pledger1971comparisons}, \parencite{li2015ordering} and \parencite{zhang2019stochastic}. In terms of the hazard rate, and the likelihood ratio orders and \parencite{amini2016skewness} established several variability orderings between order statistics from independent multiple-outlier PHR random variables. Therefore the study of stochastic ordering in the domain of Archimax copula have significant inportance for its applicability. 

In this article we will extend and prove some results on  stochastic comparision of extreme and subsequent order statistics from homogeneous random variables  in terms of generator of Archimax copula. In section two we first define some key terms related to Archimax copula and discuss how Archimax copula can give dependency structure of extreme incidents. In section three we will prove a lemma and a theorem in PHR ordering in Archimax copula. In the fourth section we will prove some stochastic ordering results for maximum and subsequent order statistics in Archimax copula. In the last we will prove some stochastic ordering results for first and second order statistics. 
\section{Preliminaries} 
\textbf{Majorization and Related orders}

Let $a_{1:n} \leq \cdots \leq a_{n:n}$ be sorted in the increasing order of the components of real vector $\boldsymbol{a}=\left(a_{1}, \ldots, a_{n}\right)$.
\subsection*{Definition $2.1.$} Suppose that  $ \boldsymbol{b}=\left(b_1, \ldots, b_n\right) \in \mathbb{R}^{n} $ is a vector. Then  the vector $\boldsymbol{b}$ is said called  weakly super-majorized by $\boldsymbol{a} \in \mathbb{R}^{n}$ ( that is, $\boldsymbol{b} \preceq^{\mathrm{w}} \boldsymbol{a}$ ) if $\sum_{i=1}^{k} b_{i:n} \geq \sum_{i=1}^{k} a_{i:n}$ for $k=1, 2, 3, \ldots, n$.

The vector $\boldsymbol{b}$ is called weakly sub-majorized by $\boldsymbol{a} \in \mathbb{R}^{n}$ (that is $\boldsymbol{b} \preceq_{\mathrm{w}} \boldsymbol{a}$ ) if $\sum_{i=1}^{k} b_{n-i+1:n} \leq \sum_{i=1}^{k} a_{n-i+1:n}$ for $k=1, 2, 3, \ldots, n$.

The vector $\boldsymbol{b}$ is said called majorized by $\boldsymbol{a} \in \mathbb{R}^{n}$ (that is $\boldsymbol{b} \preceq_{\mathrm{m}} \boldsymbol{a}$ ) if $\sum_{i=1}^{n} b_{i}=\sum_{i=1}^{n} a_{i}$ and $\sum_{i=1}^{k} b_{i:n} \geq \sum_{i=1}^{k} a_{i:n}$ for $k=1, 2, 3, \ldots, n-1$.

The vector $\boldsymbol{b}$ is said called $p$-smaller than $\boldsymbol{a} \in \mathbb{R}_{+}^{n}$ (that is $\boldsymbol{b} \stackrel{\mathrm{p}}{\leq} \boldsymbol{a}$ ) if $\prod_{i=1}^{k} b_{i:n} \geq \prod_{i=1}^{k} a_{i:n}$ for $k=1, 2, 3, \ldots, n$.

From the above definitions note that $\boldsymbol{b} \stackrel{\mathrm{p}}{\leq} \boldsymbol{a}$ and $\left(\ln b_{1}, \ldots, \ln b_{n}\right) \preceq^{\mathrm{w}}\left(\ln a_{1}, \ldots, \ln a_{n}\right)$ are equivalent and $\boldsymbol{b} \preceq_{\mathrm{m}} \boldsymbol{a}$ implies $\boldsymbol{b} \preceq^{\mathrm{w}} \boldsymbol{a}$ and $\boldsymbol{b} \preceq_{\mathrm{w}} \boldsymbol{a}$. On may refer \parencite{li2015ordering} for detailed theory on majorization.
\subsection*{Definition 2.2.} A real valued function $\Phi$ defined on a set $A \subseteq \mathbb{R}^n$ is called Schur-convex on $A$ if for $\boldsymbol{x, y} \in A$ and $\boldsymbol{x} \preceq_{\mathrm{m}} \boldsymbol{y}$ implies $\Phi(\boldsymbol{x}) \leq \Phi(\boldsymbol{y})$.
 \subsection*{Stochastic Ordering of two random Variables X and Y} Suppose $X$ and $Y$ are two random variables with CDFs $F_X$ and $F_Y$, the reliability functions $\bar{F}_{X}$ and $\bar{F}_{Y}$ and densities $f_X$ and $f_Y$ respectively. Then\newline
 $X$ is called smaller than $Y$ in the stochastic order (that is  $X \leq_{\text {st }} Y$) if $\bar{F}_{X}(x) \leq \bar{F}_{Y}(x)$ for all $x$;\newline
$X$ is called smaller than $Y$ in the reversed hazard rate order (that is  $\left.X \leq_{\mathrm{rh}} Y\right)$ if $\frac{ F_{Y}(x)}{F_{X}(x)}$ is a increasing function of $x$;\newline
$X$ is called smaller than $Y$ in the hazard rate order (that is  $X \leq_{h r} Y$ ) if the ratio $\frac{\bar{F}_{Y}(x)}{\bar{F}_{X}(x)}$ is a increasing function of $x$;\newline
$X$ is called smaller than $Y$ in the likelihood ratio order (that is  $X \leq_{\operatorname{lr}} Y$ ) if the ratio $\frac{f_{Y}(x)}{f_{X}(x)}$ is a increasing  function of $x$;

One may refer to \parencite{kochar2007stochastic} for further details on stochastic order and their properties.
\section*{ Archimax copulas and related key concepts.} We begin with defining some key concepts including Archimax copula.

In the following, we denote vectors in $\mathbb{R}^{n}$ by boldface letters, namely $\boldsymbol{x}=\left(x_{1}, \ldots, x_{n}\right)$ and $\|\cdot\|$ stands for the $\ell_{1}$-norm, namely $\|\boldsymbol{x}\|=\sum_{i=1}^{n}x_{i}$. 

\subsection{Copula}
\textbf{Definition}\\
Suppose that $\boldsymbol{X}=\left(X_{1}, \ldots, X_{n}\right)$ be a random vector with distribution function $F$, survival function $\bar{F}$ and univariate marginal distribution functions $F_{1}, \ldots, F_{n}$, then a function $C_X$ is called the copula of the random vector $\boldsymbol{X}$ if,
$$
C_X\left(x_{1}, \ldots, x_{n}\right)=F\left(F_{1}^{-1}\left(x_{1}\right), \ldots, F_{n}^{-1}\left(x_{n}\right)\right), \quad 0<x_{1}, \ldots, x_{n}<1 .
$$ 
In parallel, the survival copula of the random vector $\boldsymbol{X}=\left(X_{1}, \ldots, X_{n}\right)$ is denoted by $\widehat{C}_X$ and defined as,
$$
\widehat{C}_{\boldsymbol{X}}\left(x_{1}, \ldots, x_{n}\right)=\bar{F}\left(\bar{F}_{1}^{-1}\left(x_{1}\right), \ldots, \bar{F}_{n}^{-1}\left(x_{n}\right)\right), \quad 0<x_{1}, \ldots, x_{n}<1
$$
If $F$ is continuous then copula is uniquely defined,  otherwise the copula is unique  in the product space of the domains of all marginal CDF $F_{i}$ i.e, in $\prod_{i=1}^{n}F_{i}^{-1}[0, 1]$. As the copula doesn't include information about marginal distributions, it serves as a convenient means to impose a dependency structure onto predetermined marginal distributions in practical scenarios. Copulas have gained widespread acceptance and are important tool in risk management, finance, due to their numerous effective applications.
\subsection*{Archimedean copula and generator}
A function $\phi:[0, \infty) \rightarrow[0,1]$ with the properties that, it is non-increasing, n-monotone,  continuous function with $\phi(0)=1$,  $\lim _{x \rightarrow \infty} \phi(x)=0$, strictly decreasing on $\left[0, x_{\phi}\right)$ and with pseudo inverse $\psi$ where $x_{\phi}=$ $\inf \{x: \phi(x)=0\}$, is called an Archimedean generator and $C_{\phi}$ is called an  Archimedean copula if ,\begin{equation}
    C_{\phi}\left(x_{1}, \ldots, x_{n}\right)=\phi\left(\psi\left(x_{1}\right)+\cdots+\psi\left(x_{n}\right)\right), \quad \text { for }\left(x_{1}, \ldots, x_{n}\right) \in[0,1]^{n} \text{   see \parencite{li2015ordering}} \label{1}
\end{equation} 

\textbf{Example 2.1.} Recall from \parencite{nelsen2006archimedean} the 2-copulas,\\ \text{1.} $C_{1}(x,y)=\bigg[ max \bigg(x^{-\alpha}+y^{-\alpha}-1,0\bigg)\bigg]^{-\frac{1}{\alpha}};$ having generator  $\phi_{\alpha}(t)=(1+\alpha t)^{-\frac{1}{\alpha}},\,\, \alpha \in [-1,\infty)\setminus \{0\},\,\, 0\le x,y\le 1;$ known as Clayton 2-copula family.\\
\text{2.} $C_2(x,y)=e^{-\bigg([-\ln x]^\alpha +[-\ln y]^\alpha\bigg)^\frac{1}{\alpha}}$ having generator $\phi_{\alpha}(t)=e^{-t^\frac{1}{\alpha}}, \alpha \in [1, \infty)$ and  $0\le x,y \le 1.$
Known as Gumbel copula family.\\
\text{3.} $C_3(x,y)=1-\bigg[(1-x)^\alpha + (1-y)^\alpha -(1-x)^\alpha (1-y)^\alpha\bigg]^{\frac{1}{\alpha}},\, \alpha \in [1, \infty )$ and $ 0 \le x, y \leq 1. $ Known as Joe copula family.\\
For more examples one may refer to book by R. B. Nelsen, second edition [p-116]. 
\subsection{p-monotone function} A function $\phi$ is called $p$-monotone function, for $p \in \mathbb{N}$ and $p \geq 2$, if $\phi$ is differentiable on $(0, \infty)$ up to the order $p-2$ and satisfy $(-1)^{n} \phi^{(n)}(x) \geq 0$ for $x \in(0, \infty)$, $n \in\{1, \ldots, p-2\}$ and $(-1)^{p-2} \phi^{(p-2)}$ is non-increasing and convex on $(0, \infty)$. Every Archimedean generator of a $n$-copula is a $n$-monotone function.

It is clear that 2-monotone means that $\phi$ is convex, and that a $n$-monotone function is also $p$-monotone for $p \leq n$. If a $n$-copula $C_{\phi}$ is Archimedean, then n-monotonicity of its generator $\phi$ is  necessary \parencite{mcneil2009multivariate}.
\subsection{Stable Tail Dependence Function}
 Recall from \parencite{ressel2013homogeneous} that a function $\ell: \mathbb{R}_{+}^{n} \rightarrow \mathbb{R}^{+}$is said to be $n$-variate stable tail dependence function iff it satisfies the following properties:\\
(a) $\ell$ is homogeneous function of degree 1, that is, for all $c>0$ and $v_{1}, \ldots, v_{n} \in[0, \infty)$, 
\begin{equation}  
\ell\left(c v_{1}, \ldots, c v_{n}\right)=c \ell\left(v_{1}, \ldots, v_{n}\right); \label{2}
\end{equation}
(b) $\ell\left(\mathbf{f}_{1}\right)=\cdots=\ell\left(\mathbf{f}_{n}\right)=1$ where for $j \in\{1, \ldots, n\}$, $\mathbf{f}_{j}$ is a vector with components are zero except the $j$ th equal to 1.\\
(c) $\ell$ is fully $n$-max decreasing function, that is, for any $c \in \mathbb{N}, x_{1}, \ldots, x_{n}, h_{1}, \ldots, h_{n} \in[0, \infty)$ and $J \subseteq\{1, \ldots, n\}$ with $|J|=c$.

\begin{equation}
   \sum_{\iota_{1}, \ldots, \iota_{c} \in\{0,1\}}(-1)^{\iota_{1}+\cdots+\iota_{c}} \ell\left(x_{1}+\iota_{1} h_{1} \mathbf{1}_{1 \in J}, \ldots, x_{n}+\iota_{n} h_{n} \mathbf{1}_{n \in J}\right) \leq 0 \label{3}
\end{equation}
\subsection{Extreme Value Copula}

Let $X_i = (X_{i1}, \ldots, X_{id})$, $i \in \{1, \ldots, n\}$, be a sample of independent and identically distributed (iid) random vectors with joint distribution function $F$, marginal distribution functions $F_1, \ldots, F_d$, and copula $C_F$. Assume $F$ is continuous function. Consider the component-wise maxima:
\[ L_n = (L_{n,1}, \ldots, L_{n,d}), \]
where $L_{n,j} = \max_{i=1}^n X_{ij}$, denoting maximum. Let the joint and marginal distribution functions of $L_n$ are $F_n$ and $F_{n,1}, \ldots, F_{n,d}$ respectively, it follows that the copula $C_n$ of $L_n$ is given by,
\begin{equation}
 C_n(x_1, \ldots, x_d) = \bigg[C_F\left({x_1}^\frac{1}{n}, \ldots, x_d^\frac{1}{n}\right)\bigg]^n, \quad (x_1, \ldots, x_d) \in [0, 1]^d.  \label{4}  
 \end{equation}

The family of extreme-value copulas are the limits of these copulas $C_n$ as the sample size $n$ tends to infinity. For reference see \parencite{gudendorf2010extreme}.

\textbf{Definition:} A copula $C$ is called an extreme-value copula if there is a copula $C_F$ such that
\begin{equation}  
 \bigg[C_F\left({x_1}^\frac{1}{n}, \ldots, x_d^\frac{1}{n}\right)\bigg]^n \to C(x_1, \ldots, x_d) \quad \,\,as\,\, n \to \infty  \label{5}  
 \end{equation}
for all $(x_1, \ldots, x_d) \in [0, 1]^d$. The copula $C_F$ is called in the domain of attraction of $C$.
\subsection*{Example 2.4.}
Examples of extreme value copulas includes Gumbel copula, Negative logistic model or Galambos copula, Hüsler-Reiss copula and The t-EV Copula etc. For brevity we are not going into details, one may refers to \parencite{gudendorf2010extreme}.
\subsection{Archimax copula} Let $\phi$ be an Archimedean generator with generalized inverse $\psi$  and $\ell$ be a n-variate stdf defined in subsection $2.3$, then \parencite{mesiar2013d} proposed a suitable $n$-variate generalization of Archimax copula by setting, for all $x_{1}, \ldots, x_{n} \in[0,1]$,
\begin{equation}
C_{\phi, \ell}\left(x_{1}, \ldots, x_{n}\right)=\phi \circ \ell\left(\psi\left(x_{1}\right), \ldots, \psi\left(x_{n}\right)\right) \label{6}
\end{equation}

Now in particular, when $\ell (x)=x_{1}+\cdots+x_{n}$, that is, if $C_{\phi,\ell}$ is Archimedean, the n-monotonicity condition of $\phi$ is necessary \parencite{mcneil2009multivariate}. However, this condition is sufficient but not necessary for $C_{\phi,\ell}$ to be Archimax copula. Due to property (a) of $\ell,$ any stdf $\ell$ is uniquely defined by its restriction $A$ to the unit simplex $\Delta_{n}=\left\{\boldsymbol{w} \in[0,1]^{n}: w_{1}+\cdots+\right.$ $\left.w_{n}=1\right\}$, known as the Pickands dependence function \parencite{pickands1989multivariate}. Further, for any $\boldsymbol{x} \in \mathbb{R}_{+}^{n},\,\, \ell(\boldsymbol{x})=\|\boldsymbol{x}\| A(\boldsymbol{x} /\|\boldsymbol{x}\|)$. For any $\boldsymbol{x}=(x_1,x_2,\ldots,x_n) \in[0,1]^{n}$ and using the norm  $\|\boldsymbol{x}\|=\sum_{i=1}^{n}x_{i}$ we have , 

$$
C_{\phi, A}(\boldsymbol{x})=\phi[\|\psi(\boldsymbol{x})\| A\{\psi(\boldsymbol{x}) /\|\psi(\boldsymbol{x})\|\}]
$$
That is,\begin{equation}
C_{\phi, A}(x_1,x_2,\cdots,x_n)=\phi\bigg[\bigg( \sum_{i=1}^{n}\psi(x_i)\bigg )\cdot A\bigg(\frac{\psi(x_1)}{ \sum_{i=1}^{n}\psi(x_i)},\frac{\psi(x_2)}{ \sum_{i=1}^{n}\psi(x_i)},\cdots,\frac{\psi(x_n)}{ \sum_{i=1}^{n}\psi(x_i)}\bigg)\bigg] \label{7} 
\end{equation}
 Archimax copulas have a extreme-value attractor. Recall that a function $g: \mathbb{R}_{+} \rightarrow \mathbb{R}_{+}$ is regularly varying with index $\theta \in \mathbb{R}$ iff for all $x>0$, $g(x t) / g(t) \rightarrow x^{\theta}$ as $t \rightarrow \infty$, in notation $g \in \mathcal{R}_{\theta}$. When $1-\phi(1 / x) \in \mathcal{R}_{-\theta}$ for $\theta \in(0,1]$, it is shown in Proposition 6.1 of \parencite{charpentier2014multivariate} that $C_{\phi, \ell}$ is in the maximum domain of attraction of the extreme-value copula $C_{\ell_{\alpha}}$.
\subsection*{Example 2.5.}
We will consider here Gumbel-Hougaard copula, which is an Archimedean copula as well as extreme value copula. The function $\ell$ in this copula will be our stable tail dependece function and construction process is as follows, \\
Consider an Archimedian copula 
$$
C_{\phi}\left(x_{1}, \ldots, x_{n}\right)=\phi \left(\psi\left(x_{1}\right)+\cdots+\psi\left(x_{n}\right)\right), \quad\left(x_{1}, \ldots, x_{n}\right) \in[0,1]^{n}
$$
with generator $\phi:[0,\infty] \rightarrow[0, 1]$ and inverse $\psi(t)=\inf \{u \in[0,1]: \phi(u) \leqslant t\}$; the function $\phi$ must be strictly decreasing and convex and satisfy $\phi(0)=1$, and $\phi$ should be $n$-monotone on $(0, \infty)$ and the stable tail dependence function for Gumbel-Hougaard Copula being, 
$$
\ell\left(x_{1}, \ldots, x_{n}\right)= \begin{cases}\left(x_{1}^{\theta}+\cdots+x_{n}^{\theta}\right)^{1 / \theta} & \text { if } 1 \leqslant \theta<\infty \\ x_{1} \vee \cdots \vee x_{n} & \text { if } \theta=\infty\end{cases} 
$$
for $\left(x_{1}, \ldots, x_{n}\right) \in[0, \infty)^{n}$ see \parencite{genest1989characterization}. The Archimax copula associated to $\ell$ and $\phi$ are found according to equation \eqref{6} as , 
\begin{equation}
C_{\phi}\left(x_{1}, \ldots, x_{n}\right)=\phi \left [\left(\left(\psi( x_{1})\right)^{\theta}+\cdots+\left(\psi( x_{n})\right)^{\theta}\right)^{1 / \theta}\right ] \,\,\,for\,\, \theta \in [1,\infty). \label{8} 
\end{equation}
This is the required Archimax copula, where $\ell$ is same as stdf of the Gumbel-Hougaard or logistic copula. Now we will have an Archimax $n$-copula corresponding to each pairs of a $n$-monotone Archimedean generator and a stdf $\ell$ .\\
\textbf{1. } Gumbel copula with parameter $\alpha$, $\psi(t)=(-\ln t)^\alpha$ and $\phi(t)=e^{-t^{1/\alpha}}$ ; it easy to verify that $\phi$ is completely monotone $(i.e., (-1)^k\phi^{k}(t) \geq 0, \,\,k=0,1,2,\cdots)$, hence an $n-$monotone generator,  where $\alpha \in[1, \infty)$.\\
Therefore the requited Archimax copula is $$C_{\phi}\left(x_{1}, \ldots, x_{n}\right)= e^{-\left [\left(\left(-\ln{x_{1}}\right)^{\theta\alpha}+\cdots+\left(-\ln{x_{n}}\right)^{\theta\alpha}\right)^{1 / \theta\alpha}\right ]} ;\,\,\, \theta,\alpha \in [1,\infty).$$\\
\textbf{2. } From  example $4.1$ of \parencite{nelsen2006archimedean} the Archimedean generator $\phi(t)=\left(1+t\right)^{-\alpha},$ $\alpha \in(0,\infty);$ with inverse $\psi(t)=t^{-\frac{1}{\alpha}}-1$. Note that $\phi$ is an completely monotone Archimedean generator. Therefore the Archimax copula is 
$$C_{\phi}\left(x_{1}, \ldots, x_{n}\right)=\bigg(1+\left((x_{1}^{-\frac{1}{\alpha}}-1)^{\theta}+\cdots+(x_{n}^{-\frac{1}{\alpha}}-1)^{\theta}\right)^{1 / \theta}\bigg)^{-\alpha} ,\,\,\alpha \in (0, \infty),\,\, \theta \in [1, \infty).$$
\textbf{3. } The Clayton copula generator $\phi(t)=\left(1+\alpha t\right)^{-\frac{1}{\alpha}}$ with inverse $\psi(t)=\frac{1}{\alpha}\bigg(t^{-\alpha}-1\bigg)$. Here $\phi$ is completely monotone function,  hence an $n$-monotone Archimedean generator. Therefore the Archimax copula is $$C_{\phi}\left(x_{1}, \ldots, x_{n}\right)=\left(1+\alpha \left(\bigg[\frac{1}{\alpha}\bigg(x_1^{-\alpha}-1\bigg)\bigg]^{\theta}+\cdots+\bigg[\frac{1}{\alpha}\bigg(x_n^{-\alpha}-1\bigg)\bigg]^{\theta}\right)^{1 / \theta}\right)^{-\frac{1}{\alpha}}, \,\,\,\theta,\alpha \in [1,\infty).$$
\subsection{Archimedean and Extreme value copulas are particular case of Archimax copulas}
The class of Archimax copula is larger class that contains both Archimedean copula class as well as extreme-value copula class. When $\ell$ is the function $\ell(\boldsymbol{x})=x_{1}+$ $\cdots+x_{n}$ for all $\boldsymbol{x} \in \mathbb{R}_{+}^{n}$, then the Archimax copula $C_{\phi, \ell}$ in the equation \eqref{6} becomes the Archimedean copula,
\begin{equation}
C_{\phi}\left(x_{1}, \ldots, x_{n}\right)=\phi\left\{\psi\left(x_{1}\right)+\cdots+\psi\left(x_{n}\right)\right\} 
\label{9}
\end{equation}
When $\phi(t)=e^{-t}$ for any $t \geq 0$,  then $C_{\phi, \ell}$ in the equation \eqref{6} gives extreme-value copula $C_{\ell}$ with stdf $\ell$. In particular when $\ell=\ell_{L}$ that is maximum of the components. $C_{\phi, \ell_{L}}$ is the Fréchet-Hoeffding upper bound irrespective of $\phi$, this copula describes dependency structure between co-monotonic random variables, see \parencite{caperaa2000bivariate} and \parencite{charpentier2014multivariate}.
\section{Sample Maximum of random variables with Proportional Hazards} 
In this section we focus on maximum of two random vectors having proportional hazards and dependency described by Archimax copulas.

Random vector $(X_{1}, \ldots, X_{n})$ is said to follow proportional hazards model ( that is PHR $(\bar{B}, \boldsymbol{\alpha})$ ) if $X_{i}$ have survival function  $\bar{F}_X(x)=(\bar{B}(x))^{\alpha_{i}}$, where $\bar{B}$ is the baseline survival function  and $\boldsymbol{\alpha}=\left(\alpha_{1}, \ldots, \alpha_{n}\right)$ with $\alpha_{i}>0$ for $i=1, \ldots, n,$ is parameter vector. Assume that random vector have the Archimax copula with generator $\phi$ and pickands dependence function $A$. We denote the model by $\operatorname{PHR}(\bar{B}, \boldsymbol{\alpha} ; \phi,A)$ and let $X_{n: n}(\boldsymbol{\alpha} ; \phi,A)$ as $n$-th order statistics. 
 \subsection{Lemma A.1.} 
  Let $X_1,X_2,....,X_n$ be random variables, $\phi_1$ and $\phi_2$ are generators of $n$-dimentional Archimax copulas $C_1$ and $C_2$ respectively with inverse of $\phi_2$ is $\psi_2$. If $\psi_2\circ\phi_1$ satisfies the following,
  
$(1)$ Super-additive.

$(2)$ For $A$ such that,  $\frac{1}{n}\le A\le1 ,$
\begin{equation} \psi_2\circ\phi_1(t\cdot A)\ge A\psi_2\circ\phi_1(t), \label{10} \end{equation} where $ A:[0,1]^n \to [\frac{1}{n},1] $ is pickands dependence function as in equation \eqref{7}.

then $C_1(\textbf{x})\le C_2(\textbf{x})$ for all $\textbf{x} \in [0,1]^n$.
\subsection*{Proof:} For any $\boldsymbol{x} \in[0,1]^{n}$, by the super-additivity of $\psi_{2} \circ \phi_{1}$ we have, 
$$
\sum_{i=1}^{n} \psi_{2}\left(x_{i}\right)=\sum_{i=1}^{n} \psi_{2} \circ \phi_{1}\left(\psi_{1}\left(x_{i}\right)\right) \leq \psi_{2} \circ \phi_{1}\left(\sum_{i=1}^{n} \psi_{1}\left(x_{i}\right)\right)
$$
$$
\implies \sum_{i=1}^{n} \psi_{2}\left(x_{i}\right)\leq \psi_{2} \circ \phi_{1}\left(\sum_{i=1}^{n} \psi_{1}\left(x_{i}\right)\right)
$$ Now for $A\equiv A\bigg (\frac{\psi(x_1)}{\sum_{i=1}^{n}\psi(x_i)},\frac{\psi(x_2)}{\sum_{i=1}^{n}\psi(x_i)},\cdots,\frac{\psi(x_n)}{\sum_{i=1}^{n}\psi(x_i)}\bigg )\in [\frac{1}{n},1]$ we get,
$$ A \sum_{i=1}^{n} \psi_{2}\left(x_{i}\right)\leq A \psi_{2} \circ \phi_{1}\left(\sum_{i=1}^{n} \psi_{1}\left(x_{i}\right)\right)
$$
Now using the above inequality and inequality \eqref{10}  we get,
$$ A \sum_{i=1}^{n} \psi_{2}\left(x_{i}\right)\leq \psi_{2} \circ \phi_{1}\left(A \sum_{i=1}^{n} \psi_{1}\left(x_{i}\right)\right)
$$
and hence from the monotonic decreasing property of $\phi_{2}$ it follows that, 
$$
C_{2}(\boldsymbol{x})=\phi_{2}\left(A\bigg (\frac{\psi_{2}(x_1)}{\sum_{i=1}^{n}\psi_{2}(x_i)},\frac{\psi_{2}(x_2)}{\sum_{i=1}^{n}\psi_{2}(x_i)},\cdots,\frac{\psi_{2}(x_n)}{\sum_{i=1}^{n}\psi_{2}(x_i)}\bigg )\sum_{i=1}^{n} \psi_{2}\left(x_{i}\right)\right) $$ $$\geq \phi_{1}\left(A\bigg (\frac{\psi_{1}(x_1)}{\sum_{i=1}^{n}\psi_{1}(x_i)},\frac{\psi_{1}(x_2)}{\sum_{i=1}^{n}\psi_{1}(x_i)},\cdots,\frac{\psi_{1}(x_n)}{\sum_{i=1}^{n}\psi_{1}(x_i)}\bigg )\sum_{i=1}^{n} \psi_{1}\left(x_{i}\right)\right)=C_{1}(\boldsymbol{x})
$$
This completes the proof.
\subsection{Examples of generators satisfying the condition (2) in Lemma A.1}
 Consider two extreme value copulas Joe copula and Gumbell copula with generators respectively,

  $\phi_1(t)=1-(1-e^{-t})^{\frac{1}{\theta_1}}$ \, and \, $\phi_2(t)=e^{-t^\frac{1}{\theta_2}}$ for $\theta_1,\,\, \theta_2 \in [1,\infty]$. Therefore \,\,  $\psi_2(t)=(-\log t)^{\theta_2}$ \,\,and \,\, $\psi_2\circ\phi_1(t)=\left(-\log\left(1-\left(1-e^{-t}\right)^{\frac{1}{\theta_1}}\right)\right)^{\theta_2}$. 
Then from graphic software it can be shown that $\psi_2\circ\phi_1(t)$ satisfies the below property $(2)$ namely, for $A$ with,\, $\frac{1}{n}\le A\le1$, $\psi_2\circ\phi_1(tA)\ge A\psi_2\circ\phi_1(t)$ if  $\theta_2=1$ and $\theta_1>\theta_2$.

We will now present the Theorem for PHR stochastic ordering.
\subsection{Theorem 3.1.} Suppose $\left(X_{1}, \ldots, X_{n}\right)$ following $ \operatorname{PHR}\left(\bar{B}, \boldsymbol{\alpha} ; \phi_{1},A\right)$ and $\left(Y_{1}, \ldots, Y_{n}\right)$ following $ \operatorname{PHR}\left(\bar{B}, \boldsymbol{\beta} ; \phi_{2},A\right)$ where $\phi_1$ and $\phi_2$ are Archimedean $n-$monotone generators for Archimax copula $C_{\phi_1,A}$ and $C_{\phi_2,A}$ respectively and $\bar{B}$ the baseline servival function of the model. Then $X_{n: n}\left(\boldsymbol{\alpha} ; \phi_{1},A\right) \geq_{\text {st }} Y_{n: n}\left(\boldsymbol{\beta} ; \phi_{2},A\right)$ if \textbf{ (i) } $\beta \le ^w \alpha$ , $\psi_2\circ\phi_1$ is super additive and \textbf{ (ii) } $\psi_2\circ\phi_1(t\cdot A)\ge A \cdot\psi_2\circ\phi_1(t)$ .
\subsection*{Proof:} We know that $0\le \bar{B}(x)\le 1$ so, $\ln \bar{B}(x)$ is negative, $\boldsymbol{\beta} \preceq^{\mathrm{w}} \boldsymbol{\alpha}$ gives $\ln \bar{B}(x) \cdot \boldsymbol{\beta} \preceq_{\mathrm{w}} \ln \bar{B}(x) \cdot \boldsymbol{\alpha}$. Further, $1-e^{x}$ is decreasing and concave in $x \in \mathbb{R}$. By Theorem $(iv)$ of  \parencite{marshall1979inequalities} ch-5, it gives, 
$$
\left(1-\bar{B}^{\beta_{1}}(x), \ldots, 1-\bar{B}^{\beta_{n}}(x)\right) \preceq^{\mathrm{w}}\left(1-\bar{B}^{\alpha_{1}}(x), \ldots, 1-\bar{B}^{\alpha_{n}}(x)\right)
$$
 from the decreasing convexity of $\psi_{2}$ and Theorem $(iii)$ of \parencite{marshall1979inequalities} ch-5, it follows that,
$$
\left(\psi_{2}\left(1-\bar{B}^{\beta_{1}}(x)\right), \ldots, \psi_{2}\left(1-\bar{B}^{\beta_{n}}(x)\right)\right) \preceq_{\mathrm{w}}\left(\psi_{2}\left(1-\bar{B}^{\alpha_{1}}(x)\right), \ldots, \psi_{2}\left(1-\bar{B}^{\alpha_{n}}(x)\right)\right)
$$
Then, from the definition of weakly majorization above ordering implies,
$$
\sum_{i=1}^{n} \psi_{2}\left(1-\bar{B}^{\beta_{i}}(x)\right) \leq \sum_{i=1}^{n} \psi_{2}\left(1-\bar{B}^{\alpha_{i}}(x)\right)  
$$
$$\implies A \sum_{i=1}^{n} \psi_{2}\left(1-\bar{B}^{\beta_{i}}(x)\right) \leq A  \sum_{i=1}^{n} \psi_{2}\left(1-\bar{B}^{\alpha_{i}}(x)\right) $$
where $A\equiv A\bigg (\frac{\psi(u_1)}{\sum_{i=1}^{n}\psi(u_i)},\frac{\psi(u_2)}{\sum_{i=1}^{n}\psi(u_i)},\cdots,\frac{\psi(u_n)}{\sum_{i=1}^{n}\psi(u_i)}\bigg )\in [\frac{1}{n},1].$\\ Now the monotonicity of $\phi_{2}$ implies 
\begin{equation}\label{33(1)}
   \phi_{2}\left(A \sum_{i=1}^{n} \psi_{2}\left(1-\bar{B}^{\beta_{i}}(x)\right)\right) \geq \phi_{2}\left(A \sum_{i=1}^{n} \psi_{2}\left(1-\bar{B}^{\alpha_{i}}(x)\right)\right) 
\end{equation}
Therefore, using Lemma $A.1$ and the super-additivity of $\psi_{2} \circ \phi_{1}$ we get,
\begin{equation}\label{33(2)}
\phi_{2}\left(A \sum_{i=1}^{n} \psi_{2}\left(1-\bar{B}^{\alpha_{i}}(x)\right)\right) \geq \phi_{1}\left(A \sum_{i=1}^{n} \psi_{1}\left(1-\bar{B}^{\alpha_{i}}(x)\right)\right)
\end{equation}
Consequently, using \eqref{33(1)} and \eqref{33(2)} it holds that
$$
\begin{aligned}
F_{Y_{n: n}\left(\boldsymbol{\beta} ; \phi_{2},A\right)}(x) & =\phi_{2}\left(A \sum_{i=1}^{n} \psi_{2}\left(1-\bar{B}^{\beta_{i}}(x)\right)\right) \\
& \geq \phi_{2}\left(A \sum_{i=1}^{n} \psi_{2}\left(1-\bar{B}^{\alpha_{i}}(x)\right)\right) \\
& \geq \phi_{1}\left(A \sum_{i=1}^{n} \psi_{1}\left(1-\bar{B}^{\alpha_{i}}(x)\right)\right) \\
& =F_{X_{n: n}\left(\boldsymbol{\alpha} ; \phi_{1}, A\right)}(x), 
\end{aligned}
$$
Hence $X_{n: n}\left(\boldsymbol{\alpha} ; \phi_{1}, A\right) \geq_{s t} Y_{n: n}\left(\boldsymbol{\beta} ; \phi_{2}, A\right)$.

This completes the proof.
\subsection*{Corollary 3.1. } Since Archimedean copula is a particular type of Archimax copula when $A(\textbf{u})\equiv 1$, hence theorem $3.1$ for stochastic ordering of highest order staistics of two samples from PHR models is holds true. For proof replace $A(\textbf{u})\equiv 1$ and proceed as proof of theorem $3.1$ or see \parencite{li2015ordering}. 
\section*{Stochastic Ordering of Homogeneous Random Variables In Archimax Copulas}
Here we will state and prove some stochastic ordering property of order statistics from homogeneous random variables in terms of generator of Archimax copula. For results in Archimedean copula see \parencite{li2015ordering}.
\section{Sample Maximum  with an Archimax copula}
We employ Archimax copula to analyse the sample of homogeneous random variables. Subsequently, we establish the stochastic ordering properties of sample maximum with its adjacent order statistics. At first, we will present the Theorem 4.1.

\subsection*{Theorem 4.1.} Suppose homogeneous random variables $X_{1}, \ldots, X_{n+1}$ with Archimax copula $C_{\phi, A}$, $n$-monotone generator $\phi$ and Pickands dependence function $A$. Then ,

(i)  $X_{n-1: n} \leq_{\mathrm{rh}} X_{n: n} \leq_{\mathrm{rh}} X_{n+1: n+1} \iff \frac{t \phi^{\prime}(t)}{ \phi(t)}$ is decreasing function of $t \geq 0$.

(ii)  $X_{n-1: n} \leq_{\mathrm{hr}} X_{n: n} \leq_{\mathrm{hr}} X_{n+1: n+1} \iff \frac{t \phi^{\prime}(t)}{ (1-\phi(t))}$ is increasing function of $t \geq 0$.

(iii)  $X_{n-1: n} \leq_{\operatorname{lr}} X_{n: n} \leq_{\operatorname{lr}} X_{n+1: n+1} \iff \frac{t \phi^{\prime \prime}(t)}{\phi^{\prime}(t)}$ is decreasing function of $t \geq 0$. 
\subsection*{Example 4.1.}
\textbf{(i) }Consider for illustration of the theorem $4.1$ , the sample $X_1,X_2,X_3,X_4,X_5 $ from $U(0 , 1)$ with $n=4$ and $n+1=5$. Now we know that $X_{k:n} \sim Beta(k, n+1-k)$. Now we take $\ell$ from  Gumbel-Hougaard Copula in \parencite{charpentier2014multivariate} example $(2.5)$ as the appropriate stdf for Archimax copula model with 5-monotone Archimedean generator $\phi(t)=e^{-t^{\frac{1}{\theta}}}$, $\theta \in[1,\infty)$,  where stable tail dependence function (stdf) is $\ell(\textbf{x} )=(x_{1}^\theta + x_{2}^\theta + \cdots +x_{n}^\theta )^{\frac{1}{\theta}}$ with $\ell(\textbf{x})=\|\textbf{x}\|A(\textbf{x}/\|\textbf{x}\|)$. For reverse hazard rate  ordering, we calculate  $t \phi^{\prime}(t) / \phi(t) = -\frac{t^{\frac{1}{\theta}}}{\theta}$ which is a decreasing function of $t \ge 0 $, because for $\theta \in [1,\infty)$, we get $\frac{d}{dt}\bigg(- \frac{t^{\frac{1}{\theta}}}{\theta}\bigg)=-\frac{t^{\frac{1}{\theta}-1}}{\theta^2} \le 0 $ for all $t \ge 0$. On the other hand $\frac{F_{X_{5:5}}}{F_{X_{4:4}}}=x , \frac{F_{X_{(4:4)}}}{F_{X_{3:4}}}=\frac{x^4}{4x^3-3x^4} $ and equating both expression of $\ell$ we have $A=n^{\frac{1}{\theta}-1}$ with $n=4$. Plots are shown below which illustrates the theorem 4.1 (i) for $\theta=4$.
\begin{figure}[ht]
    \centering
    \includegraphics[width=1\linewidth]{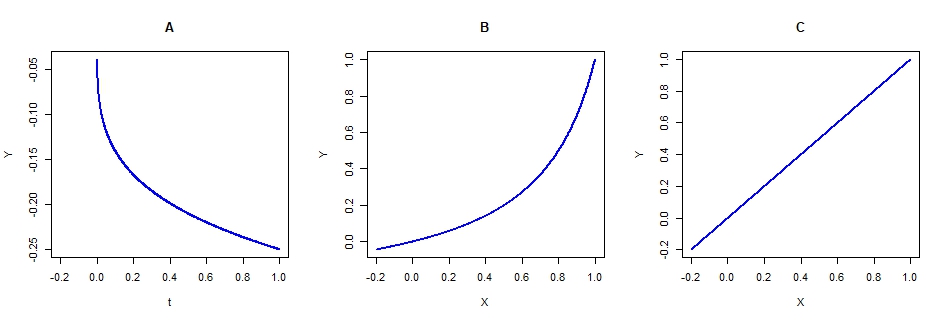} 
    \caption{Graphs of (A)\, $\frac{t \phi^{\prime}(t)}{\phi(t)} = -\frac{t^{\frac{1}{\theta}}}{\theta}, \, \theta=4, \,A=4^{\frac{1}{4}-1}$\,\,\,\,\,(B)\, $\frac{F_{X_{4:4}}}{F_{X_{3:4}}}$ \,\,\,\,\,(C)\, $\frac{F_{X_{5:5}}}{F_{X_{4:4}}}.$}
    \label{fig:1}
\end{figure} \newpage
\textbf{(ii) }For hazard rate ordering, we calculate  $t \phi^{\prime}(t) /(1- \phi(t)) = -\frac{t^{\frac{1}{\theta}}}{\theta (e^{t^{\frac{1}{\theta}}}-1)}$ which is an increasing function of $t \ge 0 $, because for  $\theta \in [1,\infty)$, we get $\frac{d}{dt}\bigg(-\frac{t^{\frac{1}{\theta}}}{\theta (e^{t^{\frac{1}{\theta}}}-1)}\bigg)=\frac{t^{\frac{1}{\theta}-1}\left(e^{t^{\frac{1}{\theta}}}t^{\frac{1}{\theta}}-e^{t^{\frac{1}{\theta}}}+1\right)}{\theta^2 (e^{t^\frac{1}{\theta}}-1)^2} \ge 0 $ for all $t \ge0$. On the other hand $\frac{1-F_{X_{5:5}}}{1-F_{X_{4:4}}}=\frac{1-x^5}{1-x^4} , \frac{1-F_{X_{4:4}}}{1-F_{X_{3:4}}}=\frac{4x^3}{1+3x^4-4x^3}$. Plots are shown below which illustrates the theorem 4.1 (ii) for $\theta=8$. 
\begin{figure}[ht] 
     \centering
     \includegraphics[width=1\linewidth]{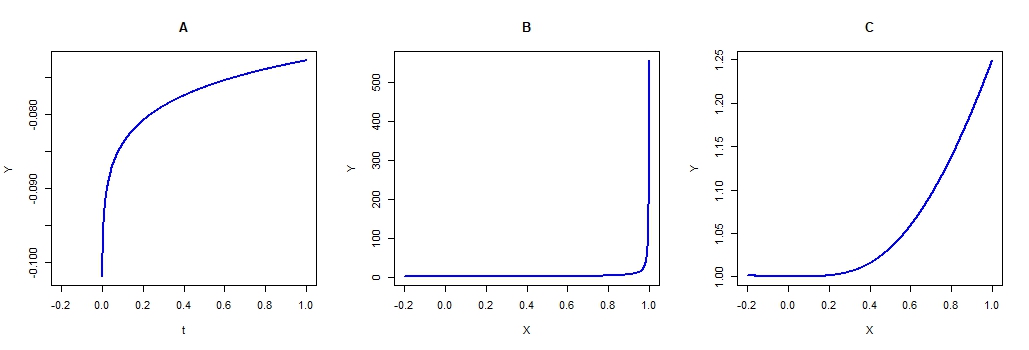}
     \caption{Graphs of (A)\, $\frac{t \phi^{\prime}(t)}{1- \phi(t)} = -\frac{t^{\frac{1}{\theta}}}{\theta (e^{t^{\frac{1}{\theta}}}-1)},\, \theta=8,\,A=4^{\frac{1}{8}-1}$ \,\,\,\,\,(B)\, $\frac{1-F_{X_{4:4}}}{1-F_{X_{3:4}}}$ \,\,\,\,\,(C)\,  $\frac{1-F_{X_{5:5}}}{1-F_{X_{4:4}}}.$ }
     \label{fig:2}
 \end{figure}

\textbf{(iii). } For likelihood ratio ordering, we calculate  $t \phi^{\prime\prime}(t) / \phi^{\prime}(t) = \frac{1-t^{\frac{1}{\theta}}}{\theta}-1$,  which is a decreasing function of $t \ge 0 $, because for $\theta \in [1,\infty)$, we get $\frac{d}{dt}\bigg(\frac{1-t^\frac{1}{\theta}}{\theta}-1\bigg)=\frac{-t^{\frac{1}{\theta}-1}}{\theta^2} \le0$ for all $t\ge0$. On the other hand $\frac{f_{X_{5:5}}}{f_{X_{4:4}}}=\frac{5x}{4} , \frac{f_{X_{4:4}}}{f_{X_{3:4}}}=\frac{4x^3}{3x^2-3x^3}$. Plots are shown below which illustrates the theorem 4.1 (iii) for $\theta=5$.
\begin{figure}[ht] 
    \centering
    \includegraphics[width=1\linewidth]{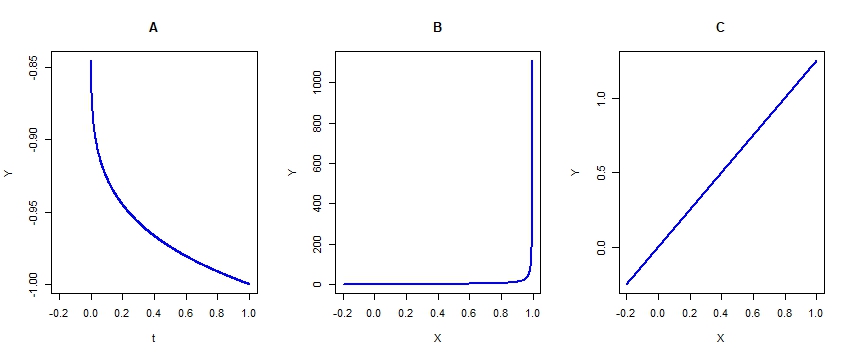}
    \caption{ Graphs of (A)\, $\frac{t \phi^{\prime\prime}(t)}{ \phi^{\prime}(t)} =\frac{1-t^{\frac{1}{\theta}}}{\theta}-1,\, \theta=5,\, A=4^{\frac{1}{5}-1}\,\,\,\,\,$(B)\,   $\frac{f_{X_{4:4}}}{f_{X_{3:4}}}$\,\,\,(C)\, $\frac{f_{X_{5:5}}}{f_{X_{4:4}}}.$ }
    \label{fig:3}  
\end{figure}

\subsection*{Proof: (i) } We first prove 
$ X_{n: n} \leq_{\mathrm{rh}} X_{n+1: n+1}$ if and only if $ \frac{t \phi^{\prime}(t)}{ \phi(t)}$ is decreasing  function of $t \geq 0$; \\
 Let $X_{i} \sim F$ for $i \in $ $\{1,2,\cdots,n+1\}$ and the Archimax Copula $C_{\phi, A}$ defined in equation \eqref{8} by,
 $$C_{\phi, A}(u_1,u_2,\cdots,u_{n+1})=\phi\bigg[\bigg( \sum_{i=1}^{n+1}\psi(u_i)\bigg )\cdot A\bigg(\frac{\psi(u_1)}{ \sum_{i=1}^{n+1}\psi(u_i)},\frac{\psi(u_2)}{ \sum_{i=1}^{n+1}\psi(u_i)},\cdots,\frac{\psi(u_{n+1})}{ \sum_{i=1}^{n+1}\psi(u_i)}\bigg)\bigg]; $$
 When we have homogeneous random variables and $u_1=u_2= \cdots = u_n = F$ then for brevity we write, $$A\bigg(\frac{\psi(F)}{ \sum_{i=1}^{n+1}\psi(F)},\frac{\psi(F)}{ \sum_{i=1}^{n+1}\psi(F)},\cdots,\frac{\psi(F)}{ \sum_{i=1}^{n+1}\psi(F)}\bigg)=A\bigg(\frac{\psi(F)}{ (n+1)\psi(F)},\frac{\psi(F)}{ (n+1)\psi(F)},\cdots,\frac{\psi(F)}{ (n+1)\psi(F)}\bigg)$$ $$=
 A\bigg(\frac{1}{ (n+1)},\frac{1}{ (n+1)},,\cdots,\frac{1}{ (n+1)}\bigg)=A_{n+1},\,\,which\,\, depends\,\,on\,\,n.$$ And we use the notation $u=\psi(F(x)) \geq 0$ and the variable $x$ in $F(x)$ is suppressed for brevity.
 Then distribution functions of $X_{n: n}$ and $X_{n+1: n+1}$ are respectively, 
$$
\begin{aligned}
& F_{X_{n: n}}(x)=\mathrm{P}\left(X_{1} \leq x, \ldots, X_{n} \leq x\right)=C_{\phi, A}(u,u,\cdots,u)=\phi[n \psi(F)A_{n}],  \\
& F_{X_{n+1: n+1}}(x)=\mathrm{P}\left(X_{1} \leq x, \ldots, X_{n+1} \leq x\right)=C_{\phi, A}(u,u,\cdots,u)=\phi[(n+1) \psi(F)A_{n+1}] ;
\end{aligned}
$$ 
Now consider,
\begin{align}  &\bigg[\frac{F_{n+1: n+1}(x)}{F_{n: n}(x)}\bigg]'  =\left[\frac{\phi[(n+1)u\cdot A_{n+1}]}{\phi[n u\cdot A_n]}\right]^{\prime} \label{11}\\
& \stackrel{\text{sign}}{=} \left[\frac{\phi^{\prime}[(n+1)uA_{n+1}](n+1) uA_{n+1}}{\phi[(n+1) uA_{n+1}]}-\frac{\phi^{\prime}[nu A_n] n uA_n}{\phi[n uA_n]}\right] \cdot \frac{\psi^{\prime}(F) f(x) \phi[(n+1) uA_{n+1}]}{u \phi[n u]} \nonumber \\ 
& \stackrel{\text{sign}}{=}\, \frac{n uA_n \phi^{\prime}[n uA_n]}{\phi[n uA_n]}-\frac{(n+1) uA_{n+1} \phi^{\prime}[(n+1) uA_{n+1}]}{\phi[(n+1) uA_{n+1}]} \nonumber
\end{align}
which is non-negative if and only if $ \frac{t\phi^{\prime}(t)}{\phi(t)}$ is decreasing function of $t\geq0$. Hence $X_{n: n} \leq_{\mathrm{rh}} X_{n+1: n+1}$ if and only if $\frac{t\phi^{\prime}(t)}{\phi(t)}$ is decreasing function of $t\geq0$.\\
Now According to $(3.4.3)$ of \parencite{david2004order}, $X_{n-1: n}$ has the distribution function in terms  of Archimax copula as ,
\begin{equation}
F_{X_{n-1: n}}(x)=n \phi[(n-1) \psi(F)A_{n-1}]-(n-1) \phi[n \psi(F)A_n];\label{12}
\end{equation}
Consider, 
$$
\begin{aligned}
\bigg[\frac{F_{X_{n-1: n}}(x)}{F_{X_{n: n}}(x)}\bigg]^\prime &=\bigg[\frac{n \phi[(n-1) uA_{n-1}]}{\phi[n uA_{n}]}-(n-1)\bigg]^\prime \\
&\stackrel{\text{sign}}{=}\bigg[ \frac{\phi[(n-1) uA_{n-1}]}{\phi[n uA_{n}]}  \bigg]'.
 \end{aligned}
 $$
Further similar simplification as after equation \eqref{11} gives, 
$$
\begin{aligned}\bigg[\frac{F_{X_{n-1: n}}(x)}{F_{X_{n: n}}(x)}\bigg]^\prime
 &\stackrel{\text{sign}}{=}\frac{n uA_n \phi^{\prime}[n uA_n]}{\phi[n uA_n]}-\frac{(n-1) uA_{n-1} \phi^{\prime}[(n-1) uA_{n-1}]}{\phi[(n-1) uA_{n-1}]} 
\end{aligned}
$$ 
 which is non-positive if and only if $\frac{t\phi^{\prime}(t)}{\phi(t)}$ is decreasing function of $t$. Hence $X_{n-1: n} \leq_{\mathrm{rh}} X_{n: n}$ if and only if $\frac{t\phi^{\prime}(t)}{\phi(t)}$ is decreasing function of $t$. Combining both parts proves the theorem $4.1$ $(i)$.

\textbf{Proof: (ii)}. Consider, 
\begin{align}
&\bigg[\frac{\bar{F}_{X_{n+1: n+1}}(x)}{\bar{F}_{X_{n: n}}(x)}\bigg]'= \left[\frac{1-\phi[(n+1) \psi(F)A_{n+1}]}{1-\phi[n \psi(F)A_n]}\right]^{\prime} \label{13} \\
&\stackrel{\text{sign}}{=} 
\bigg [\frac{\phi^{\prime}[(n+1)uA_{n+1}](n+1)uA_{n+1}}{1-\phi[(n+1) uA_{n+1}]}-\frac{\phi^{\prime}[nuA_n]nuA_n)}{(1-\phi[n uA_n])^2}\bigg ]\cdot \frac{(1-\phi[(n+1) uA_{n+1}])(-\psi^{\prime}(F))f(x)}{(1-\phi[n uA_n])\psi(F)}. \nonumber  \\ 
& \stackrel{\text{sign}}{=}
\frac{\phi^{\prime}[(n+1)uA_{n+1}](n+1)uA_{n+1}}{1-\phi[(n+1) uA_{n+1}]}-\frac{\phi^{\prime}[nuA_n]nuA_n}{1-\phi[n uA_n]} \nonumber \end{align}
which is non-negative if and only if $\frac{t\phi^{\prime}(t)}{1-\phi(t)}$  is increasing in $t$. Hence $X_{n: n} \leq_{\mathrm{hr}} X_{n+1: n+1}$ if and only if $\frac{t\phi^{\prime}(t)}{1-\phi(t)}$  is increasing in $t$. \\ 
Consider,  $$
\begin{aligned} 
 \bigg[\frac{\bar{F}_{X_{n-1: n}}(x)}{\bar{F}_{X_{n: n}}(x)}\bigg]'&= \bigg[\frac{1-n \phi((n-1) \psi(F))+(n-1) \phi(n \psi(F))}{1-\phi(n \psi(F))}\bigg]'\\ 
& \stackrel{\text{sign}}{=}\bigg[\frac{n[1-\phi[(n-1)uA_n)]}{1-\phi[nuA_n]}-(n-1)\bigg]' \\
& \stackrel{\text{sign}}{=}\bigg[\frac{1-\phi[(n-1)uA_n]}{1-\phi[nuA_n]}\bigg]'
\end{aligned}$$
Further similar simplification as after equation \eqref{13} gives, 
$$ \bigg[\frac{\bar{F}_{X_{n-1: n}}(x)}{\bar{F}_{X_{n: n}}(x)}\bigg]'\stackrel{\text{sign}}{=} \frac{\phi^{\prime}[(n-1)uA_{n-1}](n-1)uA_{n-1}}{1-\phi[(n-1) uA_{n-1}]}-\frac{\phi^{\prime}[nuA_n]nuA_n}{1-\phi[n uA_n]} $$ 
which is non-positive if and only if $\frac{t\phi^{\prime}(t)}{1-\phi(t)}$ is increasing in $t$. Hence $X_{n-1: n} \leq_{\mathrm{hr}} X_{n: n}$ if and only if $\frac{t\phi^{\prime}(t)}{1-\phi(t)}$ is increasing in $t$. Now combining both part we get the theorem $4.1(ii)$.

\textbf{Proof: (iii)}. Consider, 
$$ 
\begin{aligned}
 \left[\frac{f_{X_{n+1: n+1}}(x)}{f_{X_{n: n}}(x)}\right]^{\prime}&=\left[\frac{[ \phi((n+1) \psi(F)A_{n+1})]^{\prime}}{ [\phi(n \psi(F)A_n)]^{\prime}}\right]^{\prime}\\
 &\stackrel{\text{sign}}{=}\left[\frac{ \phi^{\prime}[(n+1) \psi(F)A_{n+1}]}{ \phi^{\prime}[n \psi(F)A_n]}\right]^{\prime} \,\,\,\text{ [ After simplification ]}\\
\end{aligned} $$
Further similar simplification as after equation \eqref{11} gives,
 $$\left[\frac{f_{X_{n+1: n+1}}(x)}{f_{X_{n: n}}(x)}\right]^{\prime}=\frac{n uA_n \phi^{\prime\prime}[n uA_n]}{\phi^{\prime}[n uA_n]}-\frac{(n+1) uA_{n+1} \phi^{\prime\prime}[(n+1) uA_{n+1}]}{\phi^{\prime}[(n+1) uA_{n+1}]}$$ 
which is non-negative if and only if $\frac{t \phi^{\prime \prime}(t)}{ \phi^{\prime}(t)}$ is decreasing in $t$. Hence $ X_{n: n} \leq_{\operatorname{lr}} X_{n+1: n+1}$ if and only if $\frac{t \phi^{\prime \prime}(t)}{ \phi^{\prime}(t)}$ is decreasing in $t$.\\
 Now consider,
$$ 
\begin{aligned}
 \left[\frac{f_{X_{n-1: n}}(x)}{f_{X_{n: n}}(x)}\right]^{\prime}  &\stackrel{\text{sign}}{=} \left[\frac{\bigg (n\phi[(n-1)\psi(F)A_{n-1}]-(n-1)\phi[n \psi(F)A_n]\bigg )^{\prime}}{\bigg (\phi[n \psi(F)A_n]\bigg )^{\prime}}\right]^{\prime} \\
 &\stackrel{\text{sign}}{=}\left[\frac{\phi^{\prime}[(n-1)uA_{n-1}]}{\phi^{\prime}[nuA_n]}\right]^{\prime} \,\, \text{[ After simplification ]}
\end{aligned}
$$
Further similar simplification as after equation \eqref{11} gives,
$$
\begin{aligned}
\left[\frac{f_{X_{n-1: n}}(x)}{f_{X_{n: n}}(x)}\right]^{\prime} &\stackrel{\text{sign}}{=} \frac{n uA_n \phi^{\prime\prime}[n uA_n]}{\phi^{\prime}[n uA_n]}-\frac{(n-1) uA_{n-1} \phi^{\prime\prime}[(n-1) uA_{n-1}]}{\phi^{\prime}[(n-1) uA_{n-1}]}   \end{aligned}
$$
Which is non-positive if and only if $\frac{t\phi^{\prime \prime}(t)}{\phi^{\prime}(t)}$ is decreasing function of $t$ and hence $X_{n-1: n} \leq_{\operatorname{lr}} X_{n: n}$ if and only if $\frac{t\phi^{\prime \prime}(t)}{\phi^{\prime}(t)}$ is decreasing function of $t$. Combining both parts we get theorem $4.1$ $(iii)$.
\subsection*{Corollary 4.1} Since every Archimedean copula and Extreme value copula are particular case of Archimax copula,  therefore stochastic ordering Theorem $4.1$ is holds true in the case of Archimedean and extreme value copula. For proof regarding archimedean copula put $A(\bar u)\equiv1$ in equation \eqref{7} or see \parencite{li2015ordering} and for proof regarding extreme value copula put $\phi(t)=e^{-t}$ in equation \eqref{7} and proceed as the proof of theorem $4.1$. 

\section{Sample minimum with an Archimax survival copula}
For the sample minimum, we employ the survival copula to describe the stochastic ordering of order statistics of homogeneous random variables. This subsection explores the minimum  of observations and its neighboring order statistics within a homogeneous sample associated with Archimax survival copulas. Similar to Theorem 4.1, we will establish various stochastic ordering among the sample minimum and its neighboring order statistics. For Archimedean servival copulas, one can refer to \parencite{arnold2008first}, \parencite{li2015ordering}.
\subsection*{Theorem 5.1} Suppose homogeneous random variables $X_{1}, \ldots, X_{n+1}$ with the Archimax survival copula $C_{\phi,A}$ , generator $\phi$ and pickands dependence function $A$. Then ,

(i) $X_{1: n+1} \leq_{\mathrm{hr}} X_{1: n} \leq_{\mathrm{hr}} X_{2: n} \iff \frac{t \phi^{\prime}(t)}{ \phi(t)}$ is decreasing function of $t \geq 0$.

(ii) $X_{1: n+1} \leq_{\mathrm{rh}} X_{1: n} \leq_{\mathrm{rh}} X_{2: n} \iff \frac{t \phi^{\prime}(t)}{(1-\phi(t))}$ is increasing function of $t \geq 0$.

(iii) $X_{1: n+1} \leq_{\operatorname{lr}} X_{1: n} \leq_{\operatorname{lr}} X_{2: n} \iff \frac{ t\phi^{\prime \prime}(t)}{ \phi^{\prime}(t)}$ is decreasing function of $t \geq 0$. 
\subsection*{Example $5.1$ } To illustrate theorem $5.1$, we will consider the random sample same as example $4.1$ of section $4$, Archimax sevival copula model with generator  $\phi(t)=(1+t)^{-\theta}$ where $\theta \in[1, \infty)$, stable tail dependent function $\ell$ defined as $\ell(\textbf{x} )=(x_{1}^\theta + x_{2}^\theta + \cdots +x_{n}^\theta )^{\frac{1}{\theta}}$, where $\ell(\textbf{x})=\|\textbf{x}\|A(\textbf{x}/\|\textbf{x}\|)$. \\
\textbf{(i) } For hazard rate ordering, we calculate  $t \phi^{\prime}(t) / \phi(t) = \frac{\theta}{1+t}-\theta$ which is a decreasing function of $t\ge0$, because for $\theta \in[1, \infty)$, we get $\frac{d}{dt}\bigg(\frac{\theta}{1+t}-\theta\bigg)=-\frac{\theta}{(1+t)^2} \le0$ for all $t\ge0$. On the other hand $\frac{1-F_{X_{1:4}}}{1-F_{X_{1:5}}}=\frac{1}{1-x} ,  \frac{1-F_{X_{2:4}}}{1-F_{X_{1:4}}}=\frac{1+8x^3-3x^4-6x^2}{(x-1)^4} $  and equating both expression of $\ell$ we have $A=n^{\frac{1}{\theta}-1}$ with $n=4$. Plots are shown below which illustrates the theorem 5.1 (i) for $\theta=4$.
 \begin{figure}[ht] 
     \centering
     \includegraphics[width=1\linewidth]{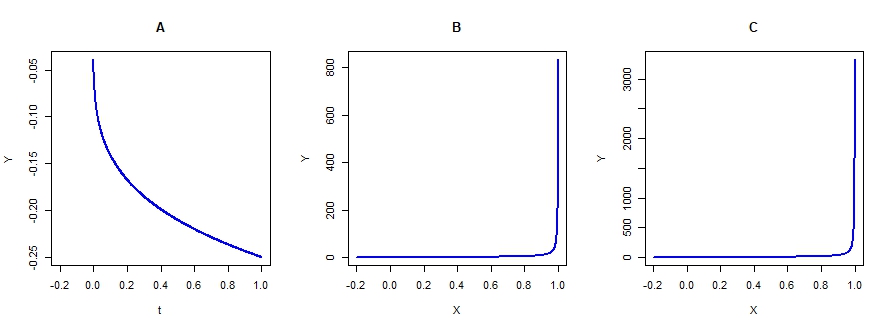}
     \caption{Graphs in  (A)\, $\frac{t\phi^{\prime}(t)}{\phi(t)} = \frac{\theta}{1+t}-\theta$,\, $\theta=4,\, A=4^{\frac{1}{4}-1}$ \,\,\,\, (B)\, $\frac{1-F_{X_{1:4}}}{1-F_{X_{1:5}}}$ \,\,\,\, (C)\, $\frac{1-F_{X_{2:4}}}{1-F_{X_{1:4}}}.$}
     \label{fig:4}
 \end{figure} \newpage
\textbf{(ii) }For reverse hazard rate ordering, we calculate  $t \phi^{\prime}(t) /(1- \phi(t)) = \frac{t\theta}{1+t-(1+t)^{\theta+1}}$ which is an increasing function of $t\ge0$, because for $\theta \in[1, \infty)$, we have  $\frac{d}{dt}\bigg(\frac{t\theta}{1+t-(1+t)^{\theta+1}}\bigg)=\frac{\theta}{1+t-\left(1+t\right)^{\left(\theta+1\right)}}-\frac{t\theta\left(1-(\theta+1)(t+1)^{\theta}\right)}{\left(1+t-\left(1+t\right)^{\left(\theta+1\right)}\right)^{2}} \ge0$ for all $t\ge0$. On the other hand $\frac{F_{X_{1:4}}}{F_{X_{1:5}}}=\frac{1-(x-1)^4}{1+(x-1)^5}$ and $ \frac{F_{X_{2:4}}}{F_{X_{1:4}}}=\frac{3x^4-8x^3+6x^2}{1-(x-1)^4} $. Plots are shown below which illustrates the theorem 5.1 (ii) for $\theta=8$.
\begin{figure}[ht] 
    \centering
    \includegraphics[width=1\linewidth]{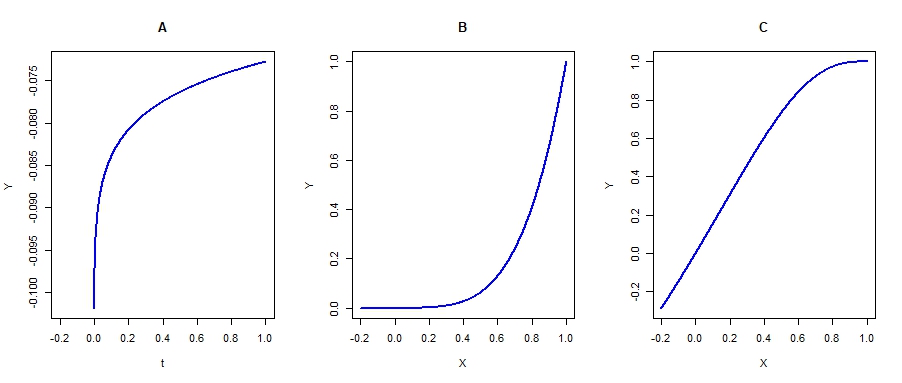}
    \caption{ Graphs in (A)\, $\frac{t \phi^{\prime}(t)}{(1- \phi(t))} = \frac{t\theta}{1+t-(1+t)^{\theta+1}},\, \theta=8,\, A=4^{\frac{1}{8}-1} $ \,\,\,\, (B)\, $\frac{F_{X_{1:4}}}{F_{X_{1:5}}}$\,\,\,\, (C)\, $ \frac{F_{X_{2:4}}}{F_{X_{1:4}}}.$}
    \label{fig:5}
\end{figure} \, \\
\textbf{(iii) } For likelihood ratio ordering, we calculate  $t \phi^{\prime\prime}(t) / \phi^{\prime}(t) = \frac{1+\theta}{1+t}-\theta-1$ which is a decreasing function of $t\ge0$, because  for $\theta \in[1, \infty)$, we have  $\frac{d}{dt}\bigg(\frac{1+\theta}{1+t}-\theta-1\bigg)=-\frac{1+\theta}{(1+t)^2}\le0$ for all $t\ge0$. On the other hand $\frac{f_{X_{1:4}}}{f_{X_{1:5}}}=\frac{4}{5(1-x)} , \frac{f_{X_{2:4}}}{f_{X_{1:4}}}=\frac{3x}{1-x}$. Plots are shown below which illustrates the theorem 5.1 (iii) for $\theta=5$. \newpage
\begin{figure}[ht] 
    \centering
    \includegraphics[width=1\linewidth]{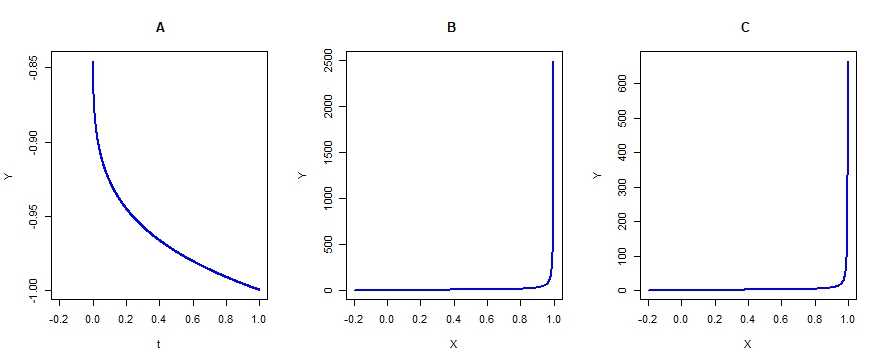}
    \caption{ Graph of (A)\, $\frac{t\phi^{\prime\prime}(t)}{\phi^{\prime}(t)} = \frac{1+\theta}{1+t}-\theta-1,\,\, \theta=5,\, A=4^{\frac{1}{5}-1},$ \,\,\,\, (B)\, $\frac{f_{X_{2:4}}}{f_{X_{1:4}}}$ \,\,\, (C)\,  $\frac{f_{X_{1:4}}}{f_{X_{1:5}}}.$}
    \label{fig:6}
\end{figure} 
\textbf{Proof:(i) } From \parencite{dos2005copula} we have,
$$
\begin{aligned} F_{X_{1:n}}(x)=1- \bar{F}_{X_{1:n}}(x)&=1-C_{\phi,A}(\bar{F}^{-1}(u),\bar{F}^{-1}(u),\cdots,\bar{F}^{-1}(u))\\
&=1-\phi[n\psi(\bar{F})A_n] \\
Similarly\,\,\, F_{X_{1:n+1}}(x)=1- \bar{F}_{X_{1:n+1}}(x)& =1-\phi[(n+1)uA_{n+1}], \,\,where\,\,u=\psi(\bar{F}(x)).
\end{aligned}
$$
Consider,
$$ 
\begin{aligned}
& \bigg[ \frac{\bar{F}_{X_{1:n}}(x)}{\bar{F}_{X_{1:n+1}}(x)}   \bigg]^{\prime} = \bigg[\frac{\phi[nuA_n]}{\phi[(n+1)uA_{n+1}]}\bigg]^{\prime} \\
&\stackrel{\text{sign}}{=}
\bigg[\frac{nA_n\psi'(\bar{F})(-f(x))\phi^{\prime}[nuA_n]}{\phi[(n+1)uA_{n+1}]}     -\frac{\phi[nuA_n]\phi^{\prime}[(n+1)uA_{n+1}](n+1)A_{n+1}\psi'(\bar{F})(-f(x))}{(\phi[(n+1)uA_{n+1}])^2} \bigg].
\end{aligned}
$$
 
$\stackrel{\text{sign}}{=} \bigg[\frac{nuA_n\phi^{\prime}[nuA_n]}{\phi[nuA_n]} -\frac{\phi^{\prime}[(n+1)uA_{n+1}](n+1)A_{n+1}}{\phi[(n+1)uA_{n+1}]}  \bigg] $  [ After simplification ]\\  
which is a non-negative if and only if $\frac{t \phi^{\prime}(t)}{\phi(t)}$ is decreasing in $t$. Hence $X_{1: n+1} \leq_{\mathrm{hr}} X_{1: n}$ if and only if $\frac{t \phi^{\prime}(t)}{\phi(t)}$ is decreasing in $t$. \\
From \parencite{david2004order},  $\bar{F}_{X_{2:n}}(x)=n\phi[(n-1)uA_{n-1}]-(n-1)\phi[nuA_n] $ \\
Now consider, 
$$
\begin{aligned}
 \bigg[ \frac{\bar{F}_{X_{2:n}}}{\bar{F}_{X_{1:n}}}   \bigg]^{\prime}& = \left[\frac{n\phi[(n-1)uA_{n-1}]-(n-1)\phi[nuA_n]}{\phi[nuA_n]}\right]^{\prime} \\
&\stackrel{\text{sign}}{=} \bigg[\frac{\phi[(n-1)uA_{n-1}]}{\phi[nuA_n]}\bigg]' \text{ [ After simplification ]}  
\end{aligned}
$$
Further similar simplification as after equation \eqref{11} gives, $$ \bigg[ \frac{\bar{F}_{X_{2:n}}}{\bar{F}_{X_{1:n}}}   \bigg]^{\prime}\stackrel{\text{sign}}{=}\frac{(n-1) uA_{n-1} \phi^{\prime}[(n-1) uA_{n-1}]}{\phi[(n-1) uA_{n-1}]}-\frac{n uA_n \phi^{\prime}[n uA_n]}{\phi[n uA_n]}$$\\ which is non-negative if and only if  $\frac{t \phi^{\prime}(t)}{ \phi(t)}$ is decreasing function of $t$. Hence $ X_{1: n} \leq_{\mathrm{hr}} X_{2: n}$ if and only if  $\frac{t \phi^{\prime}(t)}{ \phi(t)}$ is decreasing function of $t$. Combining both parts we get theorem $5.1$ $(i)$.\\
\textbf{Proof: (ii) } Consider,
    $$
    \begin{aligned}
 \bigg[\frac{F_{X_{1:n}}}{F_{X_{1:n+1}}}   \bigg]^{\prime} &\stackrel{\text{sign}}{=} \bigg[ \frac{1-\phi[nuA_n]}{1-\phi[(n+1)uA_{n+1}]}\bigg]^{\prime} \\
 &\stackrel{\text{sign}}{=}\frac{\phi^{\prime}[(n+1)uA_{n+1}](n+1)uA_{n+1}}{1-\phi[(n+1) uA_{n+1}]}-\frac{\phi^{\prime}[(n)uA_n]nuA_n}{1-\phi[nuA_n]} \,\,[ \mbox{After \ simplification} ]\end{aligned}
$$
 
which is non-negative if and only if $\frac{t\phi^{\prime}(t)}{(1-\phi(t))}$ is increasing function of $t$. Hence $X_{1: n+1} \leq_{\mathrm{rh}} X_{1:n} $ if and only if $\frac{t\phi^{\prime}(t)}{(1-\phi(t))}$ is increasing function of $t$.\\ 
Now consider,
$$
\begin{aligned}
\bigg[ \frac{F_{X_{2:n}}}{F_{X_{1:n}}}   \bigg]^{\prime} &= \left[\frac{1-n\phi[(n-1)uA_{n-1}]+(n-1)\phi[nuA_n]}{1-\phi[nuA_n]}\right]^{\prime} \\
& \stackrel{\text{sign}}{=} \bigg[\frac{1-\phi[(n-1)uA_{n-1}]}{1-\phi[nuA_n]}\bigg]'\,\,\text{[ After simplification ] }\\
&\stackrel{\text{sign}}{=}\frac{\phi^{\prime}[nuA_n]nuA_n}{1-\phi[n uA_n]}-\frac{\phi^{\prime}[(n-1)uA_{n-1}](n-1)uA_{n-1}}{1-\phi[(n-1) uA_{n-1}]} \ge 0.\,\, [\, Further\,\,\, simplification\,]\end{aligned} $$
which is non-negative if and only if $\frac{t\phi^{\prime}(t)}{(1-\phi(t))}$ is increasing function of $t$. Hence $ X_{1: n} \leq_{\mathrm{rh}} X_{2: n}$ if and only if $\frac{t\phi^{\prime}(t)}{(1-\phi(t))}$ is increasing function of $t$.
Hence combining both the results we get theorem $5.1(ii)$ .

\textbf{Proof: (iii) } Now consider,
$$ 
\begin{aligned}
 \left[\frac{f_{X_{1: n}}(x)}{f_{X_{1: n+1}}(x)}\right]^{\prime}&=
\left[\frac{[1- \phi(n \psi(\bar{F})A_n)]^{\prime}}{ [1- \phi((n+1) \psi(\bar{F})A_{n+1})]^{\prime}}\right]^{\prime} \\ 
 &\stackrel{\text{sign}}{=}  \left[\frac{ \phi^{\prime}[nuA_n]}{ \phi^{\prime}[(n+1) uA_{n+1}]}\right]^{\prime}
\end{aligned}$$
Further similar simplification as after equation \eqref{11} gives,

$$\left[\frac{f_{X_{1: n}}(x)}{f_{X_{1: n+1}}(x)}\right]^{\prime}\stackrel{\text{sign}}{=}\frac{n uA_n \phi^{\prime\prime}[n uA_n]}{\phi^{\prime}[n uA_n]}-\frac{(n+1) uA_{n+1} \phi^{\prime\prime}[(n+1) uA_{n+1}]}{\phi^{\prime}[(n+1) uA_{n+1}]} $$

which is non-negative if and only if $\frac{t \phi^{\prime \prime}(t)}{ \phi^{\prime}(t)}$ is decreasing function of $t$ and hence $X_{1: n+1} \leq_{\operatorname{lr}} X_{1: n}$ if and only if $\frac{t \phi^{\prime \prime}(t)}{ \phi^{\prime}(t)}$ is decreasing function of $t$.\\
Now consider,
$$
\begin{aligned}
\bigg[ \frac{f_{X_{2:n}}}{f_{X_{1:n}}}   \bigg]^{\prime}
 &= \left[\frac{\bigg(1-n\phi[(n-1)uA_{n-1}]+(n-1)\phi[nuA_n]\bigg)'}{\bigg(1-\phi[nuA_n]\bigg)'}\right]^{\prime}\\
&\stackrel{\text{sign}}{=} \bigg[\frac{\phi'[(n-1)uA_{n-1}]}{\phi'[nuA_n]}\bigg]'\cdot\frac{(n-1)A_{n-1}}{nA_n}\\
&\stackrel{\text{sign}}{=} \bigg[\frac{\phi'[(n-1)uA_{n-1}]}{\phi'[nuA_n]}\bigg]' \,\, \text{[ After simplification ]}
\end{aligned}
$$
Further similar simplification as after equation \eqref{11} gives, $$\bigg[ \frac{f_{X_{2:n}}}{f_{X_{1:n}}}   \bigg]^{\prime}\stackrel{\text{sign}}{=} \frac{(n-1) uA_{n-1} \phi^{\prime\prime}[(n-1) uA_{n-1}]}{\phi^{\prime}[(n-1) u A_{n-1}]}-\frac{n uA_n \phi^{\prime\prime}[n uA_n]}{\phi^{\prime}[n uA_n]}$$  
which is non-negative if and only if $\frac{t \phi^{\prime \prime}(t)}{  \phi^{\prime}(t)}$ is decreasing function of $t$. Hence $X_{1: n} \leq_{\mathrm{lr}} X_{2: n}$ if and only if $\frac{t \phi^{\prime \prime}(t)}{\phi^{\prime}(t)}$ is decreasing function of $t$. Combining both results we get theorem $5.1(iii)$.
\subsection*{Corollary 5.1} Since every Archimedean copula and extreme value copula are particular case of Archimax copula,  therefore stochastic ordering Theorem $5.1$ is true in the case of Archimedean and extreme value survival copula. For proof regarding Archimedean survival copula put $A(\bar u)\equiv1$ in equation \eqref{7} (one may also refer \parencite{li2015ordering}). And for proof regarding extreme value survival copula put $\phi(t)=e^{-t}$ in equation \eqref{7} and proceed as the proof of the theorem $5.1$.

\section{Some Applications}
Stochastic ordering properties of order statistics has many applications in multivariate statistics, one important aspect is k out of n system in reliability theory. An article by \parencite{bergmann1991stochastic} discussed on stochastic ordering of lifetimes of systems with independent components and \parencite{li2015ordering}
focused on the same in Archimedean copula. In this article we accomplished stochastic ordering of systems lifetime when exchangeable components lifetime dependency is described by Archimax copula. Suppose that stochastic representation of Archimax copula $C_A$ as in \parencite{charpentier2014multivariate} 
is $( X_1, X_2, \cdots , X_n )= ( RS_1, RS_2, \cdots, RS_n )$, where $R(>0)$ is radial component of the process with CDF F and $(S_1, S_2, \cdots, S_n)$ is the simplex component of the process and independent of $R$. Then we know that the $n$ monotone generator $\phi$ of the Archimax copula is given by the formula, 
$$\phi(x)=\int_{x}^{\infty}\bigg(1-\frac{x}{r}\bigg)^{n-1}dF(r)$$
Now recall that $X\leq_{\mathrm{hr}}Y$ if for any time $t$, the hazard functions satisfies $hr_Y(t)\le hr_X(t)$, that is, hazard rate of the system having lifetime $Y$ is less than that of system with lifetime $X$, provided that the systems $X$ and $Y$ survived till the time $t$. Similarly $X\leq_{\mathrm{rh}}Y$ if for any time  $t$, the reverse hazard functions satisfies $rh_X(t)\leq rh_Y(t)$, that is, reverse hazard rate of the system having lifetime $X$ is less than that of system with lifetime $Y$, provided that the systems $X$ and $Y$ survived till the time $t$.
 
Now theorem $4.1(i)$ states that $X_{n-1: n} \leq_{\mathrm{rh}} X_{n: n} \leq_{\mathrm{rh}} X_{n+1: n+1}$ that is, at any arbitrary time $t$, reverse hazard rate of 2 out of n system $\leq$ that of $n$ parallel system $\leq$ that of $(n+1)$ parallel system  if and only if $\frac{t \phi^{\prime}(t)}{ \phi(t)}$ is decreasing function of $t$, where $\phi$ is the Archimax $n$ monotone generator.

Theorem $4.1(ii)$ states that $X_{n-1: n} \leq_{\mathrm{hr}} X_{n: n} \leq_{\mathrm{hr}} X_{n+1: n+1}$ that is, at any arbitrary time $t$, hazard rate of 2 out of n system $\geq$ that of $n$ parallel system $\geq$ that of $n+1$ parallel system if and only if $\frac{t \phi^{\prime}(t)}{1- \phi(t)}$ is increasing function of $t$. 

Theorem $5.1(i)$ states that $X_{1: n+1} \leq_{\mathrm{hr}} X_{1: n} \leq_{\mathrm{hr}} X_{2: n}$, that is,  at an arbitrary time $t$, hazard rate of $n+1$ series system $\geq$ that of $n$ series system $\geq$ that of fail safe system with $n$ components if and only if $\frac{t \phi^{\prime}(t)}{ \phi(t)}$ is decreasing function of $t$.

Theorem $5.1(ii)$ states that $X_{1: n+1} \leq_{\mathrm{rh}} X_{1: n} \leq_{\mathrm{rh}} X_{2: n}$, that is, at an arbitrary time $t$, reverse hazard rate of $n+1$ series system $\leq$ that of $n$ series system $\leq$ that of fail safe system with $n$ components if and only if $\frac{t \phi^{\prime}(t)}{1- \phi(t)}$ is increasing function of $t$.

Suppose $\left(X_{1}, \ldots, X_{n}\right)$ following $ \operatorname{PHR}\left(\bar{B}, \boldsymbol{\alpha} ; \phi_{1},A\right)$ and $\left(Y_{1}, \ldots, Y_{n}\right)$ following $ \operatorname{PHR}\left(\bar{B}, \boldsymbol{\beta} ; \phi_{2},A\right)$. Now under the assumptions of Theorem $3.1$, $X_{n: n}\left(\boldsymbol{\alpha} ; \phi_{1},A\right) \geq_{\text {st }}Y_{n:n}\left(\boldsymbol{\beta} ; \phi_{2},A\right)$, that is, reliability of parallel system with components $\left(X_{1}, \ldots, X_{n}\right)$ is greater than that of parallel system with components $\left(Y_{1}, \ldots, Y_{n}\right)$ at any time $t$. \\
\begin{figure}[tp]
  \centering

  \begin{subfigure}{1\linewidth}
    \centering
    \includegraphics[width=\linewidth]{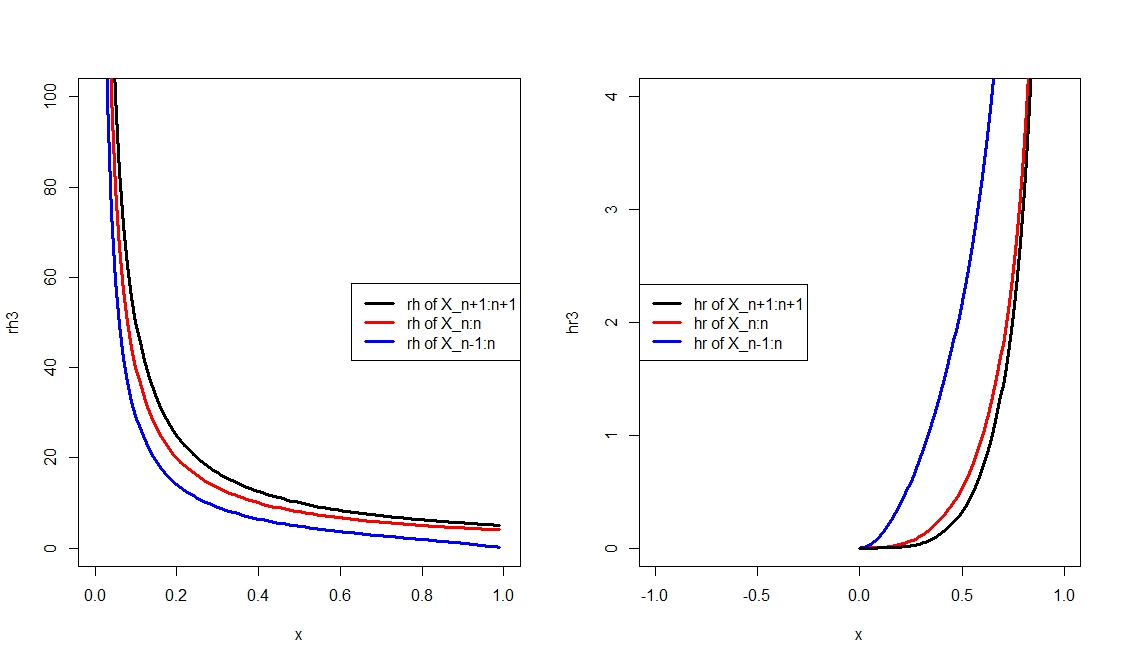}
    \caption{Graphical illustration of comparison of reverse hazard functions (Left) and hazard functions (Right) of $X_{n-1:n},X_{n:n}$ and $X_{n+1:n+1}$ for example 4.1 with $n=4$.}
    \label{fig:7}
  \end{subfigure}

  \vspace{\baselineskip} % Add some vertical space between the subfigures

  \begin{subfigure}{1\linewidth}
    \centering
    \includegraphics[width=\linewidth]{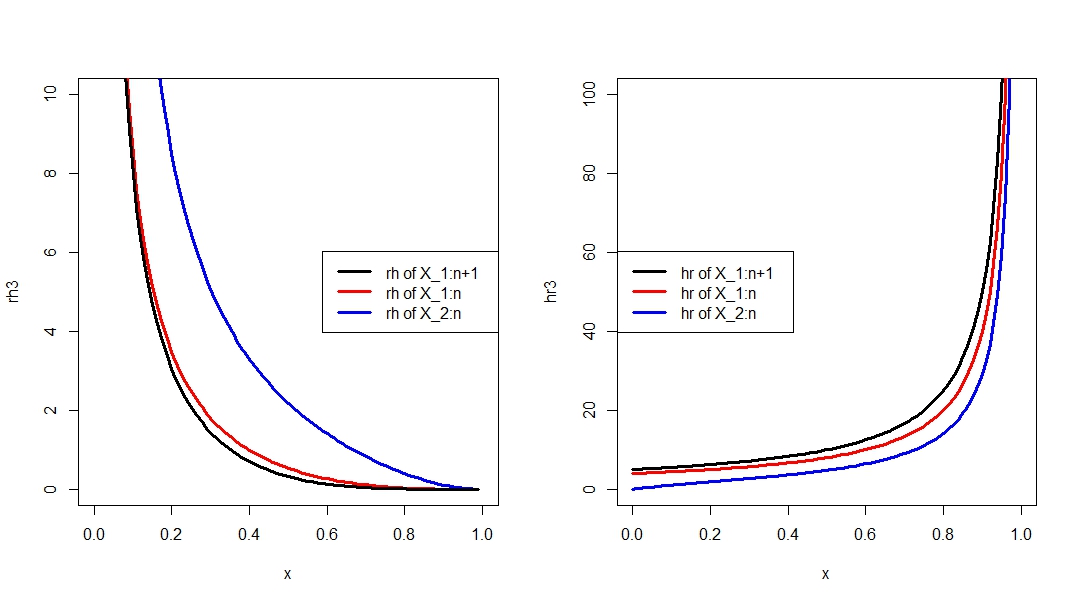}
    \caption{Graphical illustration of comparison of reverse hazard functions (Left) and hazard functions (Right) of $X_{2:n},X_{1:n}$ and $X_{1:n+1}$ for example 5.1 with $n=4$.}
    \label{fig:8}
  \end{subfigure}

  \caption{ Comparision of reverse hazard and hazard functions of extreme order statistics.}
\end{figure} 
\newpage
 \section{Conclusions} 
In this article we studied stochastic ordering of extreme order statistics when the dependency structure among the homogeneous random variables are described by Archimax copula. We also studied the stochastic ordering of maximum of two random vectors following PHR model in Archimax copula. The future of this research may lie on Stochastic ordering of order statistics from heterogeneous random variable in Archimax copula.
\section*{Acknowledgements}
Sarikul Islam is supported with research grant provided by Ministry of Human Resource Development (MHRD), Government of India. 
\printbibliography
\end{document}